\documentclass[reqno, oneside]{amsart}
\usepackage{amsmath, amssymb}

\theoremstyle{plain}
\newtheorem{Thm}{Theorem}
\newtheorem{Prop}[Thm]{Proposition}
\newtheorem{Lem}[Thm]{Lemma}
\newtheorem{Cor}[Thm]{Corollary}

\theoremstyle{definition}

\theoremstyle{remark}
\newtheorem*{Rem}{Remark}

\numberwithin{Thm}{section}

\def\fl(#1){\eqref{fl:#1}}
\def\thm#1{Theorem~\ref{thm:#1}}
\def\lem#1{Lemma~\ref{lem:#1}}
\def\cor#1{Corollary~\ref{cor:#1}}
\def\prop#1{Proposition~\ref{prop:#1}}
\def\sec#1{\S{\rm \ref{sec:#1}}}

\newcommand{\R}{\frak R}
\newcommand{\Lg}{\frak g}

\newcommand{\Sl}{\frak sl}

\newcommand{\J}{\frak J}
\newcommand{\Ug}{U_{q}(\frak g)}

\newcommand{\Uq}{U_{q}^{res}(\frak g)}

\newcommand{\Uf}{{U_q^{fin}}(\frak g)}

\newcommand{\Ud}{U_q^{fin}(\frak g)^*}

\newcommand{\Z}{\frak Z}
\newcommand{\tZ}{\widetilde{\frak Z}}
\newcommand{\pZ}{\frak Z'}

\newcommand{\Wl}{\hat{W}_l}
\newcommand{\F}{\mathcal{F}}
\newcommand{\id}{{\textrm id}}

\newcommand{\bCf}{\bar{{\mathcal C}_f}}

\newcommand{\tmX}{\bar{\mathcal X}}

\begin{document}
\author{Anna Lachowska}
\thanks{This research was supported by the NSF VIGRE grant and the Clay Mathematics Institute}
\address{Department of Mathematics\\
MIT\\
Cambridge MA 02139\\
USA}
\email{lachowska@math.mit.edu}

\title{On the Center of the Small Quantum Group}

\date\today


\begin{abstract} 

Using the quantum Fourier transform $\F$ \cite{LM}, we describe the block 
decomposition and multiplicative structure of a subalgebra $\tZ + \F(\tZ)$
 of the center of the small quantum group $\Uf$ at a root of unity. 
It contains the known subalgebra $\tZ$ \cite{BG}, which is isomorphic to
the algebra of characters of finite dimensional $\Uf$-modules. We prove 
that the intersection $\tZ \cap \F(\tZ)$ coincides with the annihilator of 
 the radical of $\tZ$. 
 Applying representation-theoretical methods, we show that 
 $\tZ$ surjects onto 
the  algebra of endomorphisms of certain indecomposable projective 
modules over $\Uf$. 
In particular this leads to the conclusion that the center of $\Uf$ 
coincides with  
 $\tZ + \F(\tZ)$ in the case $\Lg =sl_2$.

\end{abstract}

\maketitle

\section{Introduction and Notation}

Let $\Uf$ denote the small quantum group \cite{Lus2} 
associated to a semisimple complex Lie algebra $\Lg$ 
and an odd number $l$. We will assume $l$ is greater 
than the Coxeter number of the root system and relatively  
prime to the determinant of the Cartan matrix of $\Lg$.   

In the present paper we combine Hopf-algebraic and 
representation-theoretical methods to obtain the block decomposition 
and multiplicative structure of a subalgebra in the center  
$\Z$ of $\Uf$. 

The starting point is the result in \cite{BG}, which describes a 
subalgebra $\tZ \subset \Z$ (\thm{BG}).  
It can be obtained as the image of the complexification of the 
Grothendieck ring $\R$ of the category of 
finite dimensional $\Uf$-modules under the homomorphism $\J$ (\sec{R}), 
which was introduced in \cite{Dr} for any quasitriangular Hopf algebra 
and  can be viewed as an inverse of the quantum Harish-Chandra map. If 
$q$ is not a root of unity, the image of $\J$ coincides with the 
whole center of the quantum enveloping algebra $\Ug$ (see e.g. \cite{Bau}). 
 In the case of the small quantum group, $\R$ surjects onto a proper 
subalgebra of the 
center $\J(\R) \simeq \tZ$, as was already observed in the $\Sl_2$ 
case in \cite{Ker1}.
 
On the other hand, viewing $\Uf$ as a finite dimensional Hopf 
algebra allows one to define an injective map $\phi^{-1} : \R \to \Z$ 
(\sec{Hopf}), and obtain an ideal in  the center $\pZ = \phi^{-1}(\R)$. 
 A combination of the two theories yields a map 
$\F \simeq \J \circ \phi$ with the properties of a Fourier transform, 
which is a slight modification of a map introduced in \cite{LM} as 
the basic involution in a projective action of the modular 
group on $\Uf$. In \sec{inters} we use the involutive property of 
$\F$ to prove the main structural result (\thm{intersec}):  
$$\tZ \cap \pZ \simeq {\rm Ann}\,({\rm Rad} \;\tZ).$$ 
Then we show that $\tZ + \pZ = \tZ + \F(\tZ)$ is a subalgebra of the
 center, and describe its multiplicative structure and block 
decomposition (\thm{main}). It contains many elements which cannot 
be expressed in terms of any version of a Harish-Chandra map, including  
the two-sided integral of $\Uf$.

We would like to mention that the ideal $\tZ \cap \pZ$ is an 
interesting object in itself. It has many recognizable properties: its
 basis is parametrized by the orbits of a certain group of affine
 reflections in $P$, it is invariant under a suitable version of the 
Fourier transform, 
and it appears naturally in the representation theory of a quantum group at 
a root of unity. 
This is reminiscent of the algebra of characters of the fusion category  
over $\Uq$ (the Verlinde algebra) \cite{AP}, and it turns out that 
the analogy between the two structures goes even further. 
We study this structure in \cite{La}, and show that it can be considered a 
natural counterpart of the Verlinde algebra for the small quantum group. 

In \sec{act} we study the action of the center on projective $\Uf$-modules 
and prove that the subalgebra $\tZ$ surjects onto the direct sum of 
algebras of endomorphisms of certain projective modules over $\Uf$. This 
result may 
be considered a quantum root of unity version of the category 
${\mathcal O}$ 
statement \cite{Soe1}, where the center of the universal enveloping 
algebra $U(\Lg)$ surjects onto the algebra of endomorphisms of the 
projective 
cover of a simple $U(\Lg)$-module with anti-dominant highest weight. 
The surjection 
$\tZ \to End_{\Uf}P(\lambda)$ holds for all the projective modules 
$P(\lambda)$ over $U_q^{fin}(\Sl_2)$, 
which allows us to prove that $\Z = \tZ + \pZ$ in the 
$\Sl_2$ case.
We discuss the cases $\Lg = \Sl_2$ and $\Lg = \Sl_3$ in \sec{ex}. 

The constructed subalgebra is the smallest subspace of the center 
containing the obvious part $\tZ$ and invariant under the quantum Fourier 
transform; it coincides with the whole center 
in the $\Sl_2$ case. We believe that its description will serve as a 
basepoint for understanding the structure of the center $\Z$ of 
$\Uf$ in general.  

\subsection{Definitions and notation}

To simplify the arguments we will assume that the root system 
$R$ of $\Lg$ is simply laced, though this restriction is not crucial. 
Let $R_+$ denote the positive roots, 
$\Pi \subset R_+$ the simple roots,
 $Q$ the root lattice, $P$ the weight lattice, and $r = {\rm rank} \Lg$.
Fix an odd number $l \geq h$, where $h$ is 
the Coxeter number of the root system $R$. Starting from \sec{R} we assume 
that $(l, {\rm det}a_{ij})=1$ for the Cartan 
matrix $a_{ij}$ of $R$ (see \lem{GCD}).

Let $W$ be the Weyl group associated to the root system $R$. Its action on 
$P$ is given by $s_\alpha(\mu) = \mu - \langle \mu, \alpha \rangle \alpha$ 
for any 
$\alpha \in R$, $\mu \in P$,
where the form  
$\langle \cdot, \cdot \rangle$  is normalized so that $\langle \alpha, \alpha 
\rangle =2$ for $\alpha \in R$. We will also use the shifted action 
$w \cdot \mu = w(\mu +\rho) -\rho$ for any $w \in W$, $\mu \in P$, where 
$\rho = \frac{1}{2} \sum_{\alpha \in R_+} \alpha$.  
Let $\Wl$ denote the affine Weyl group generated by the reflections 
$s_{\alpha, k} \cdot \mu = s_{\alpha} \cdot \mu + kl\alpha$ for any 
$\alpha \in R_+$,
$k \in {\bf Z}$ and $\mu \in P$. 
The natural 
and shifted actions of $\Wl$ are defined by 
$s_{\alpha, k} ( \mu )= s_{\alpha} ( \mu) + kl\alpha$ and 
$s_{\alpha, k} \cdot \mu = s_{\alpha} \cdot \mu + kl\alpha$ for any 
$\mu \in P$.
The sets  $ \hat{X}= \{ \lambda \in P_+ : \langle \lambda , \alpha \rangle \leq  l \;\;{\rm for}\;{\rm all}\; \alpha \in R_+ \} $
and
$ \bar{X}= \{ \lambda \in P : 0\leq \langle \lambda + \rho, \alpha \rangle \leq  l \;\;{\rm for}\;{\rm all}\; \alpha \in R_+ \} $ 
are the fundamental domains for the non-shifted and shifted affine Weyl group actions respectively. 

Let $\Ug$ be the Drinfeld-Jimbo quantum enveloping algebra associated to a semisimple simply-laced Lie algebra $\Lg$ over ${\mathbb Q}(v)$, where $v$ 
is a formal variable. It is 
generated by the elements $\{ E_i, F_i, K_i ^{\pm 1}\}_{i=1}^r$,  with 
the standard set of 
relations (see e.g. \cite{Lus1}). Denote by $U_v(\Lg)^{\mathbb Z}$ 
the ${\mathbb Z}[v,v^{-1}]$ subalgebra of $\Ug$ generated by the divided 
powers $E_i^{(k)} = \frac{E_i^k}{[k]!}$, $F_i^{(k)} = \frac{F_i^k}{[k]!}$  
and $K_i, K_i^{-1}$ for $i =1 \ldots r$,
$k \in {\mathbb N}, k \geq 1$, where
 $[k]! = \prod_{s=1}^k \frac{v^s -v^{-s}}{v -v^{-1}}.$
Then the restricted quantum group $\Uq$ is defined by 
$\Uq = U_v(\Lg)^{\mathbb Z} \otimes_{{\mathbb Z}[v, v^{-1}]} {\mathbb Q}(q)$, 
where $ v$ maps to $q \in {\mathbb C}$, which is set to be a primitive 
$l$-th root of unity. 
By \cite{Lus1}, we have for $1 \leq i \leq r$: $E_i^l =0$, $F_i^l =0$, 
$K_i^{2l} =1$, and $K_i^l$ is central in $\Uq$. 

Let  $\Uf$  be the finite dimensional subquotient of $\Uq$ \cite{Lus2}, generated by 
$\{ E_i, F_i, K_i^{\pm 1} \}_{i=1}^r$ over ${\mathbb Q}(q)$, factorized 
over the two-sided ideal $\langle \{K_i^l -1\}_{i=1}^r \rangle $. Then  
$\Uf$  is a Hopf 
algebra of dimension $l^{{\rm dim}\Lg}$ over ${\mathbb Q}(q)$. 
We will use the same notation for the $\mathbb C$-algebras 
$\Uq$ and $\Uf$, where the 
field ${\mathbb Q}(q)$ is extended to $\mathbb C$.

\subsection*{Acknowledgements.} I would like to thank my advisor 
I. Frenkel for his continuing help and encouragement. I am very grateful 
to H.H.~Andersen, I. Gordon and W. Soergel for valuable 
discussions and suggestions, and to the referee of the first version of 
this text for many helpful remarks.  

\section{Hopf structure and the isomorphism of modules.} \label{sec:Hopf}

The properties of  $\Uf$  as a finite dimensional Hopf algebra were studied 
in e.g. \cite{Lyu}, \cite{Ker1}. 
In particular, it is unimodular, i.e. there exists a two-sided integral  
$\Lambda \in$  $\Uf$  such that $a \Lambda = \Lambda a = \varepsilon(a) 
\Lambda$, for any $a \in  \Uf$, where $\varepsilon$ is the counit for 
$\Uf$, and $\Lambda$ is 
unique up to a scalar multiple. The Hopf dual  $\Ud$  of  $\Uf$  is not 
unimodular; however the spaces of left and right integrals 
$\lambda_l, \lambda_r 
\in \Ud$  such that $p \lambda_l = p(1) \lambda_l$ and 
$\lambda_r p = p(1) \lambda_r$ for all $p \in \Ud$  are one-dimensional. 
For any $\mu = \sum_i n_i \alpha_i \in Q$ we set 
$K_\mu = \prod_i K_i^{n_i}$.  
The square of the antipode $S^2$ is an inner automorphism of  $\Uf$  
with $S^2(a) = K_{-2\rho} a K_{2 \rho}$ for any $a \in \Uf$, where 
$K_{2\rho} = \prod_{\alpha \in R_+}K_\alpha$. 

To formulate the next result we will need to turn  $\Uf$  and  $\Ud$  into 
left and right modules over themselves. Naturally  $\Uf$  is a left module 
over itself via multiplication. Define a left $\Uf$-module structure on  
$\Ud$  by
$$ a \rightharpoondown p(b) = p(S(a) b),$$ 
for any $ a,b, \in \Uf$, $p \in \Ud$. 
$\Ud$ is a right module over itself via right multiplication.  $\Uf$  can be  
given a right $\Ud$- module structure by setting 
$$ a \leftharpoonup p = \sum p(a_{(1)}) a_{(2)},$$
for any $a \in \Uf$, $p \in  \Ud$, where we have used Sweedler's 
notation
for comultiplication $\Delta(a) = \sum a_{(1)} \otimes a_{(2)}$.   
 
The next theorem is a special case of a general result valid for 
finite dimensional Hopf algebras with an antipode.

\begin{Thm} \label{thm:Hopf}
\cite{Swe}, \S 5.1; \cite{Rad}.
There exists a linear map $\phi$  from  $\Uf$  to  $\Ud$,  which is both an 
isomorphism of left $\Uf$-modules and right $\Ud$-modules, satisfying the 
following conditions:  
\begin{enumerate}
\item{$\phi(1) = \lambda_r$ is a right integral for  $\Ud$  and 
$\phi^{-1} (\varepsilon) = \Lambda$ is an integral for $\Uf$;} 
 \item{$\phi(a) = a \rightharpoondown \lambda_r$ for any $a \in  \Uf$, and 
$\phi^{-1}(p) = \Lambda \leftharpoonup p $ for any $p \in  \Ud$.
Thus  $\Ud$  is a free $\Uf$-module with basis $\lambda_r$, and  $\Uf$  is 
a free $\Ud$-module with basis $\Lambda$ with respect to the actions 
defined above.} 
\end{enumerate}
With these conditions, $\phi$ is unique up to a scalar multiple.
\end{Thm}

\begin{Cor} \label{cor:ZC} The image of the center $\Z$ of $\Uf$ under 
$\phi$ coincides with 
$C_r $, the set of all functionals $p \in \Ud$  such that 
$ p(ab) =p(b S^{-2}(a))$. In particular, ${\rm dim}(\Z)= {\rm dim}(C_r)$. 
\end{Cor}

\begin{proof}
 The right integral $\lambda_r$ is itself an element of $ C_r $ 
 ( \cite{Rad}, Thm.3).
Since $S(\Z) = \Z$, the result of the left action of $\Z$ on $\lambda_r$ is 
again in $ C_r. $ Conversly, suppose 
$a \rightharpoondown \lambda_r \in C_r$ Then 
$$\lambda_r(S(a)bc) = \lambda_r(c S^{-2}(S(a)b))= 
\lambda_r( c S^{-1}(a) S^{-2}(b) ), $$
which has to be equal to $ \lambda_r(S(a)cS^{-2}(b))= \lambda_r(cS^{-2}(b)
S^{-1}(a))$ for all $b,c \in \Uf$, and  
this holds if and only if $a \in \Z$. 
\end{proof}

\begin{Rem} The set $C_r$ has a simple interpretation in terms of 
cocommutative elements of $\Ud$. Let $C_0$ be the set of all cocommutative 
elements of $\Ud$, i.e. $\mu(ab) = \mu(ba)$ for $\mu \in C_0$ and any
$a,b \in \Uf$. Consider the functional $p \in \Ud$  defined by 
$p = K_{-2\rho} \rightharpoondown \mu$.  Using the identity $S^2 = Ad(K_{-2
\rho})$, it is easy to check that $p \in C_r$: 
$$p(ab)= \mu(abK_{2\rho})=\mu(bK_{2\rho}a) =
\mu(b K_{2\rho} a K_{-2\rho} K_{2\rho}) = p(b S^{-2}(a)).$$ 
Similarly, for a functional $t \in C_l$ defined by the condition $t(ab) =
t(b S^2(a))$, we have $t = K_{2\rho} \rightharpoondown \mu$ for some 
$\mu \in C_0$. Therefore, the sets $C_r$, $C_l$ and $C_0$ differ only 
by the action of a semisimple grouplike element of $\Uf$, and since $\Delta
(K_{2\rho})= K_{2\rho} \otimes K_{2\rho}$, they are isomorphic as 
commutative 
algebras (\cite{Dr}, \S 3) and their dimensions coincide.   
\end{Rem}

Let $\R$ be the subspace of $C_0$ spanned by the traces of simple $\Uf$-
modules. 
By \cite{Lus2} they are parametrized by the $l$-restricted weights 
$$(P/lP)_+ \equiv \{ \lambda \in P : 0 \leq \langle \lambda, \alpha_i 
\rangle <l, \alpha_i \in \Pi \},$$  so the 
dimension of $\R$ equals $l^r$. $\R$ has an algebra structure induced 
from the tensor product in the category of finite dimensional  $\Uf$-
modules. It is a commutative subalgebra in $\Ud$.
The integral form of $\R$ is the Grothendieck ring of the category of finite 
dimensional 
$\Uf$-modules.
By the above remark we can consider the isomorphic algebras of $q$- and 
$q^{-1}$-traces
of $\Uf$- modules, defined as the usual traces of $\Uf$-modules, 
shifted by the action of $K_{2\rho}$ and  $K_{-2\rho}$ respectively.
They will be denoted by $\R_l \subset C_l$ and 
$\R_r \subset C_r$. We have $\R \simeq \R_r \simeq \R_l$ as 
commutative algebras. 

Let $\pZ \subset \Z$ denote the image of $\R_r$ under the map $\phi^{-1}$.

\begin{Prop} \label{prop:pZ}
$$ \pZ \simeq {\rm Ann}({\rm Rad}\;{\Z}).$$ 
\end{Prop}

\begin{proof}
 A central 
element $z$ acts on a simple $\Uf$-module as multiplication by a scalar, 
which is nonzero whenever $z$ is not nilpotent. Therefore, as a $\Z$-
module,
$\R_r$ coincides with the maximal semisimple submodule of $C_r$, which is 
annihilated by the radical of $\Z$. 

By restriction, $\phi$ is an isomorphism of $\Z$-modules, and \thm{Hopf}
provides an isomorphism of $\Z$-module structures on $C_r$ and $\Z$, 
the latter given by multiplication. 
In particular, $\pZ= \phi^{-1}(\R_r)$ coincides with the 
annihilator of the radical (i.e. the socle) of $\Z$. 
\end{proof}

\section{Quasitriangular structure in $\Uf$.} \label{sec:R}

By \cite{Lyu},  $\Uf$  is quasitriangular with the canonical element $R$ 
satisfying 
$$ R \Delta(x) R^{-1} = \Delta^{op}(x) $$ 
for any $x \in \Uf$  and 
$$ (\Delta \otimes {\rm id}) R = R_{13} R_{23}, \;\; ({\rm id}\otimes 
\Delta) R = R_{13} R_{12},$$
where $\Delta^{op}(x) = \sum x_{(2)} \otimes x_{(1)} $, and the lower 
indices of $R$ indicate the position of $R_{(1)}$ and $R_{(2)}$.

Let $U^0$ denote the subalgebra of $\Uf$  generated by 
$\{K_i^{\pm 1}\}_{i=1 \ldots r}$. Define a bimultiplicative 
symmetric form on  $U^0$ by $\pi(K_\mu, K_\nu) = q^{ (\mu| \nu )}$, where  
$( \mu | \nu ) : Q \times Q \to {\mathbb C}$ is the bilinear pairing 
such that $(\alpha_i | \alpha_j) = a_{ij}$ for any $\alpha_i, \alpha_j \in 
\Pi$ and the Cartan matrix $a_{ij}$ of $\Lg$. 

\begin{Lem} (cf. App. B, \cite{Lyu}) \label{lem:GCD}
Suppose that for the Cartan matrix $a_{ij}$ of $\Lg$, an odd number $l$ 
satisfies the condition $(l, {\rm det}a_{ij})=1$. Then 
$\pi $ is nondegenerate on $U^0$.
\end{Lem}
\begin{proof}
 Let $K_h \in Ann \pi$, which means that $\pi(K_h, K_\alpha)=1$ for any 
$\alpha \in \Pi$. Then in general, $Ann \pi = \{K_\lambda, \lambda \in 
Q \cap lP \}$ \cite{Lyu}. 
Write $h = \sum_{i =1}^r b_i \alpha_i$ for some coefficients $b_i \in 
\mathbb Z$, and 
$$\pi(K_h, K_{\alpha_i}) =
q^{( \sum b_j \alpha_j | \alpha_i )} = q^{\sum a_{ij}b_j}=1.$$
Then $\sum a_{ij}b_j \equiv 0 ({\rm mod}l)$, for any $i =1, \ldots , r$. 
Solving the system, we get $b_j = \frac{1}{{\rm det}a_{ij}} \sum A_{ij} k_i 
l$, where $A_{ij}$ are the minors of $a_{ij}$, and $k_i \in \mathbb Z$.  
By assumption this implies $b_j \equiv 0 ({\rm mod}l)$, therefore 
$h \in lQ$ and $Q \cap lP = lQ$. 
Clearly all $K_\alpha^l \in Ann \pi$, and hence $Ann \pi $ 
is generated by 
$\{K_\alpha^l \}_{\alpha \in \Pi} $, and $\pi$ is nondegenerate 
on $U^0 = \{ K_\mu, \mu \in Q/lQ \}$.
\end{proof}

From now on we will assume that $l$ is an odd integer greater than or 
equal to the Coxeter number of the root system of $\Lg$, satisfying the 
conditions of 
\lem{GCD}. 
Then $\Uf$ is a quasitriangular Hopf algebra with the canonical 
element given by the formula
\cite{Lyu}:

$$ R = \frac{1}{|Q/lQ|} \sum_{\mu, \nu \in Q/lQ} \pi(K_\mu, K_\nu)^{-1}K_
\mu \otimes K_\nu \cdot \sum_{m_\alpha =0}^{l-1} \prod_\alpha 
\frac{q^{m_\alpha(m_\alpha -1)/2} (q-q^{-1})^{m_\alpha}}{[m_\alpha]!}
F_\alpha^{m_\alpha} \otimes E_\alpha^{m_\alpha}, $$
where $\alpha$ runs over all positive roots in a fixed order corresponding 
to the decomposition of the longest element of the Weyl group 
$w_0 \in W$.
 
$\Uf$ is also factorizable, meaning that the map 
$p \rightarrow m(p \otimes {\rm id})(R_{21}R_{12})$ 
from  $\Ud$  to  $\Uf$  is surjective. Here $m$ denotes the multiplication 
in $\Uf$. By \cite{Lyu}, a small quantum group 
is factorizable if and only if for $\Gamma = Q/lQ$, $2\Gamma =\Gamma$, 
where $2\Gamma = \{2x, x \in \Gamma \}$ 
is the group of squares. This always holds for  $\Uf$  since $l$ is odd. 

\begin{Thm} \cite{Dr}  The map $\J: \Ud  \rightarrow \Uf$  defined by 
$$\J(p) = m(p \otimes {\rm id})(R_{21}R_{12}),$$
 restricted to $C_r$, is an 
algebra homomorphism between $C_r$ and $\Z$. 
\end{Thm}

\begin{Cor} In the case of a factorizable finite dimensional Hopf algebra, the 
map 
$$\J : C_r \to \Z $$ 
is an isomorphism of commutative algebras. 
\end{Cor}
\begin{proof} 
Factorizability for a finite dimensional Hopf algebra means that 
$Ker (\J) \subset \Ud$ is trivial. Since $\J(C_r) \subset \Z$ is 
a homomorphic image of $C_r$ and 
$dim(\Z)=dim(C_r)$ by \cor{ZC}, we get an isomorphism. 
\end{proof}

Denote the image of $\R_r \subset C_r$ under $\J$ by $\tZ \subset \Z$.
By the Remark in \sec{Hopf}, we have $\R_r \simeq \R$, the algebra 
 of characters of finite dimensional 
$\Uf$- modules with multiplication induced by the tensor product. 
The algebraic structure of $\R$ can be described using the representation 
theory of 
finite dimensional modules over $\Uq$. 

Let $\bCf$ denote the category of finite dimensional modules over $\Uf$. 
For a module $M$ in $\bCf$ we will use the symbol $chM$ for the formal 
character 
$chM = \sum_{\eta \in P} ({\rm dim} M_\eta) e^\eta \in {\mathbb C}[P]^W$, 
where we write 
$e^\eta$ for the basis element in ${\mathbb C}[P]$ corresponding to 
$\eta \in P$. Here ${\mathbb C}[P]^W$ denotes the subspace of exponential 
invariants of the $W$-action in ${\mathbb C}[P]$.

\begin{Prop} 
$$  \R \simeq  {\mathbb C}[P]^W \otimes_{{\mathbb C}[lP]^W} {\mathbb C}, $$
where $ {\mathbb C}[lP]^W \longrightarrow {\mathbb C}$ is evaluation at 
$1$.
\end{Prop}

\begin{proof}
The finite dimensional simple modules of type $\bf 1$ 
(i.e. all $\{K_i^l\}_{i=1\ldots r}$ act on them by $1$) over $\Uq$ are 
parametrized by the dominant integral weights $\{L(\mu) \}_{\mu \in P_+}$. 
The algebra of characters, spanned 
by $\{{\rm ch} L(\mu) = \sum_{\nu \in P} {\rm dim}( L(\mu))_\nu e^\nu \}_
{\mu \in P_+}$, 
is isomorphic to ${\mathbb C}[P]^W$.  By Lusztig's tensor product 
theorem \cite{Lus1},
$$ L(\mu_0 + l\mu_1) \simeq L(\mu_0) \otimes L(l\mu_1), $$
where $\mu = \mu_0 + l\mu_1 \in P_+$ and $0\leq \langle \mu_0, \alpha_i 
\rangle \leq l-1$ for all 
$\alpha_i \in \Pi$. 
 By restriction, simple $\Uq$-modules of type $\bf 1$ are  
modules over $\Uf$. The subquotient  $\Uf$  of $\Uq$ acts trivially on the 
simple modules 
$L(l\mu_1)$, whose weights are all in $lP$ \cite{APW2}. Therefore, the 
algebra of characters 
$\R$ of  $\Uf$-modules  can be obtained 
from that of  $\Uq$  by setting the characters of simple modules with $l$-
multiple highest weights equal to their dimensions. This leads to the 
formula above. 
\end{proof}

By the above, the subalgebra $\tZ \equiv \J(\R_r)\subset \Z$ 
is isomorphic to $\R_r$. 

\section{The block decomposition of $\tZ$.} \label{sec:tZ}

The block decomposition of $\tZ$ can be obtained using Theorem 4.5 of 
\cite{BG}. 
To formulate the result, we need to define the action of the  
Weyl group on the restricted weight lattice 
$(P/lP)_+$.

Note that the shifted action of the affine Weyl group $\Wl$ preserves the 
$(lP-\rho)$-lattice in 
$P$:
$$ s_{\alpha, k} \cdot (l\mu -\rho) = s_\alpha(l\mu) + k l \alpha -\rho \in 
(lP-\rho).$$ 
Therefore, the shifted action of $\Wl$ on $(P/lP)_+$ can be defined by 
$$ w \bullet \lambda = w \cdot \lambda \,{\rm mod}(lP)$$
for any $\lambda \in (P/lP)_+$, $w \in \Wl$.
By definition, the result of this action belongs again to $(P/lP)_+$. 
The affine Weyl group $\Wl$ is isomorphic to $W\ltimes lQ$, where the 
elements of $lQ$ generate the translations. Since $lQ \subset lP$, it is 
enough to consider the $\bullet$-action of just the elements $w \in W$ of 
the finite Weyl group. 

To find the fundamental domain of the $W\bullet$ action, consider the 
affine group $W_l^P = W \ltimes lP$, where the translations are generated 
by the $lP$ lattice. Let $\Omega$ be the subgroup of $W_l^P$ stabilizing 
the alcove $\bar{X}$ with respect to the shifted action. 
Recall that with the assumption on the odd number $l$ given by \lem{GCD}, 
we have $lP \cap Q = lQ$. 
Then by \cite{Hum}, \S 4.5, $\Wl \cap \Omega =1$, so in fact $W_l^P$ is the 
semidirect product of $\Wl$ and $\Omega$, and $\Omega \simeq W_l^P/\Wl$. 
The order of $\Omega$ is equal to $|W_l^P/\Wl|=|P/Q|$, the index of 
connection. With the above assumption on $l$, the order of the open alcove 
$|X|$ is divisible by $|P/Q|$. 

Denote by $\tmX$ the set of representatives of the orbits of the shifted 
action of $\Omega$ in the alcove $\bar{X}$. Then 
$\tmX$ is the set of representatives of the orbits of the 
$W\bullet$ action in $(P/lP)_+$. 
Let $\mathcal X \subset \tmX$ denote the subset of regular weights in 
$\tmX$, such that their stabilizer in $W$ with respect to the $W \bullet$ 
action is trivial. 
Then $|\mathcal X| = |X|/|P/Q|$. 

Comparing to the results in \cite{APW2}, one can deduce that 
the set $\tmX$ enumerates the blocks of the category of finite dimensional 
modules over $\Uf$ (see \cite{La}). 

Similarly one can define the natural action of $\Wl$ on $(P/lP)_+$ by  
$$ w \circ \lambda = w(\lambda) \,{\rm mod}(lP)$$ 
for any $\lambda \in (P/lP)_{+}$, $w \in \Wl$. It is easy to see that 
$W\bullet$ and $W\circ$ have the same orbit structure in $(P/lP)_+$.


We will also need the restriction of the action of $\Wl$ to the finite root 
lattice $Q/lQ$. The natural (non-shifted) action of $\Wl$ on $Q \subset P$ 
can be defined by restriction from $P$, since the action preserves the root 
lattice:
$$ s_{\alpha, k} (\beta) = s_\alpha(\beta) + kl\alpha $$
for any $\alpha, \beta \in Q$, $k \in \mathbb Z$.
Then the $\circ$-action of $\Wl$ on $Q/lQ$ is defined by 
$$ w \circ \beta = w(\beta) \,{\rm mod}lQ $$
for any $\beta \in Q/lQ$. As before, it is enough to consider the $\circ$-
action on $Q/lQ$ of elements of the finite Weyl group $W$.  

Now we can state the following corollary of Theorem 4.5 in \cite{BG}. 

\begin{Thm} \label{thm:BG}
Let $\tmX$ denote the set of representatives of the orbits of the 
$W\bullet$-action in 
$(P/lP)_+$, and for any $\mu \in \tmX$ let $W_\mu$ be the stabilizer 
subgroup of $\mu$ in $W$ for this action.  
 Then 
$$ \R \simeq \oplus_{\mu \in \tmX} {\mathbb C}[P]^{W_\mu} 
\otimes_{{\mathbb C}[P]^W} {\mathbb C}, $$
where ${\mathbb C}[P]^W \longrightarrow {\mathbb C}$ is evaluation at 
$1$.
\end{Thm}

\begin{proof}
Since $R$ is assumed to be simply laced, we will always identify 
the root lattice $Q$ with $\check{Q}$, the ${\mathbb Z}$-span of the coroots 
$\check{\alpha},  \alpha \in R$. 
Denote by $H^*$ the real vector space spanned by the simple roots, and 
let $H$ be its dual space with the basis $\{h_i \}_{i=1}^r$ such that 
$\alpha_i(h)=\kappa(h_i, h)$ for any $h \in H$. Here $\kappa(\cdot, \cdot)$ 
is the restriction of the Killing form in $\Lg$. Then $Q^*$, the 
${\mathbb Z}$-span of $\{h_i \}_{i=1}^r$, is a lattice in $H$ dual to $Q$.  
 Let $T = \{ e^{2\pi i h}, h \in H{\rm mod}Q^* \}$ denote the maximal torus
 of the simply connected group $G$ with $Lie(G) =\Lg$. Then the elements
 of ${\mathbb C}[P]$ can be considered as characters on  $T$ by setting  
$e^\mu (e^{2\pi i h}) = e^{2 \pi i \mu(h)} \in {\mathbb C}$,
for any $t= e^{2 \pi i h} \in T$. 

Recall that $\R \simeq {\mathbb C}[P]^W \otimes_{{\mathbb C}[lP]^W} {\mathbb C}$. 
Then we are in the setting of Theorem 4.5 in \cite{BG} 
(in our case $\chi_s =1$ in the notations of \cite{BG}), 
and we have 
$$ \R \simeq \oplus_{\{ t \in T : t^l =1 \}/W} {\mathbb C}[P]^{W(t)} 
\otimes_{{\mathbb C}[P]^W} {\mathbb C}, $$
where the map ${\mathbb C}[P]^W  \to {\mathbb C}$ is given by the 
evaluation at $1$. 
The notation $\{ t \in T : t^l =1 \}/W$ stands for the set of 
representatives of the $W$ orbits in $\{ t \in T : t^l =1 \}$, 
which has a natural $W$ action.  
 The subgroup 
$W(t) \subset W$ is generated by the reflections $\{s_
\alpha , \alpha \in R \; | \, e^\alpha(t)=1\}$.

We can rewrite the result parametrizing the blocks by the $W \circ$-orbits 
in the restricted root lattice.  Indeed, the set $ \{t \in T : t^l =1 \}$ 
contains all $t= e^{2\pi i h}$ where $lh \in Q^*$. Therefore, it is 
parametrized by  
$h \in Q^*/lQ^*$, which under duality corresponds to $\beta \in Q/lQ$, 
with the orbits of the $W$-action in $\{ t \in T : t^l =1 \}$ translated 
to the orbits of the $W\circ$-action in $Q/lQ$.  
Thus the maximal 
ideals of $\R$ can be parametrized by $(Q/lQ)/W\circ$, the representatives
of $W\circ$-orbits in the restricted root lattice. The stabilizer subgroup 
$W(t)$ is generated by those $s_\alpha, \alpha \in R$ such that 
$e^\alpha(t) = e^{2\pi i \alpha(h)} =1$, or $\alpha(h) \in {\mathbb Z}$ 
for $h \in Q^*/lQ^*$. In the dual picture this corresponds to 
those $s_\alpha, \alpha \in R$, such that 
$\langle \alpha, \beta \rangle \equiv 0(mod l)$ for 
$\beta \in Q/lQ$.  
These are exactly the reflections $s_\alpha \in W$ that stabilize the 
element $\beta \in Q/lQ$ with respect to the $W\circ$-action. Denote 
by $W_\beta \subset W$ the subgroup generated by such reflections 
for a fixed $\beta \in (Q/lQ)/W\circ$.  
Then we get 
$$ \R \simeq \otimes_{\beta \in (Q/lQ)/W\circ} {\mathbb C}[P]^{W_\beta} 
\otimes_{{\mathbb C}[P]^W} {\mathbb C}. $$

However, if we want to consider the subalgebra $\tZ \subset \Z$, which is 
isomorphic to $\R$, 
then it will be more natural to parametrize its blocks by weights in $\tmX$. 
This agrees with the decomposition of $\Uf$ as a module over itself, where 
the block corresponding to an 
element $\mu \in \tmX$ contains only the composition factors of highest 
weights in the $W \bullet$-orbit of $\mu$.  The set $(P/lP)_+$ under 
the $W\bullet
$-action has the same orbit structure as $Q/lQ$ under the $W\circ$-action. 
Indeed, with the conditions on $l$ formulated in \lem{GCD}, 
for any $\beta \in Q/lQ$ 
there exists a unique $\lambda \in (P/lP)_+$ such that 
$\langle \beta, \alpha \rangle = \langle \lambda,  
\alpha \rangle {\rm mod}(l)$, $\alpha \in \Pi$. This provides a 
one-to-one correspondence 
between the orbits of the $W\circ$ actions in $Q/lQ$ and $(P/lP)_+$. 
Clearly the $W\bullet$-action on 
$(P/lP)_+$ has the same orbit structure. Therefore we can translate the 
parametrization by 
$(Q/lQ)/W\circ$ to $(P/lP)/W\bullet \simeq \tmX$. 
In this interpretation, the subgroup $W_\mu$ is generated by the 
reflections in $W$ 
(not necessarily simple) which stabilize the weight $\mu \in \tmX$ with 
respect to the $W\bullet$-action in $(P/lP)_+$, and we obtain the statement 
of the theorem.
\end{proof}

\begin{Rem}
Using the isomorphism $\J : \R_r \to \tZ$ and \thm{BG}, we can write $\tZ 
\simeq \oplus_{\mu \in \tmX} {\tZ}_\mu$, where ${\tZ}_\mu \equiv {\mathbb 
C}[P]^{W_\mu}\otimes_{{\mathbb C}[P]^W} {\mathbb C}$ for any 
$\mu \in \tmX$. Each block ${\tZ}_\mu$ is a local algebra by \S 4.3 in 
\cite{BG}. It is also 
Frobenius, since $\oplus_{\mu \in \tmX} {\tZ}_\mu$ can be realized as the 
center of 
a matrix algebra $\oplus_{\mu \in \tmX} {\rm Mat}_{l^{|R_+|}}({\tZ}_\mu)$ 
\cite{BG}. 
Therefore, for each $\mu \in \tmX$, 
the socle of ${\tZ}_\mu$ is one-dimensional. 
\end{Rem}

\begin{Prop} \label{prop:idemp}
The elements of $\tZ$ distinguish the blocks of $\Uf$. In other words, 
 the system $\{ 1_\mu \in \tZ_\mu \}_{\mu \in \tmX} $ of idempotents in the 
algebra $\tZ$ 
is a complete set of central simple idempotents in $\Uf$, and 
$\sum_{\mu \in \tmX} 1_\mu = 1 \in$ $\Uf$. 
\end{Prop}

\begin{proof} 
The algebra of functionals $C_r$ is mapped homomorphically onto 
$\Z$ by the isomorphism $\J$. A functional in $C_r \setminus \R_r$ is not a 
trace, and hence is zero when evaluated on any element $K_\mu \in U^0$. Any 
such functional is mapped to a nilpotent central element by $\J$. 
Therefore, since $\J$ is a bijection, all idempotents in $\Z$ have their 
pre-images in $\R_r$, or equivalently all central idempotents belong to the 
subalgebra $\tZ \subset \Z$.

As a module over itself by left multiplication, $\Uf$ decomposes into a 
direct sum of 
two-sided ideals $\Uf \simeq \oplus_{\mu \in \tmX} A_\mu$, corresponding to 
the blocks 
of the category of its finite dimensional modules. This 
 induces a decomposition of the center 
$\Z = \oplus_{\mu \in \tmX} \Z_\mu$. 
By \thm{BG}, $\tZ = \oplus_{\mu \in \tmX} \tZ_\mu$ with each $\tZ_\mu 
\subset \Z_\mu$ being 
a local algebra with the idempotent $1_\mu$. Since all idempotents are in 
$\tZ$, the element $1_\mu$ is the idempotent in 
$\Z_\mu$, and $\oplus_{\mu \in \tmX} 1_\mu = 1 \in \Z$. In particular, 
each $\Z_\mu$ is a local algebra, ${\Z}_\mu = 1_\mu + {\rm Rad}(\Z_\mu)$.
\end{proof}

\section{Intersection and sum of two central subalgebras.} 
\label{sec:inters}
Since $\pZ$ is an ideal of $\Z$, the sum $\tZ + \pZ$ is a subalgebra in the 
center. 
To describe its multiplicative structure, we need to 
 find the intersection $\pZ \cap \tZ$.
This requires a combination of the two mappings used to define $\pZ$ and 
$\tZ$. 

Define the bijective mappings $\F : \Z \longrightarrow \Z$ and 
$\F' : C_r \longrightarrow C_r$ by  
$$ \F(a) \equiv \J \circ \phi(a) = 
m((\lambda_r \leftharpoonup S(a) )\otimes {\id})(R_{21}R_{12})), $$
$$ \F'(f) \equiv \phi \circ \J(f) = 
\lambda_r \leftharpoonup S(m((f\otimes {\id})(R_{21}R_{12}))). $$

For the involutive properties discussed below, $\F$ and $\F'$ will be 
called the quantum 
Fourier transforms. This differs slightly from the original definition in 
\cite{LM}, where the Fourier transforms are defined on $\Uf$ and $\Ud$
respectively, but we will need here only their restrictions to $\Z$ and 
$C_r$. 

\begin{Thm} \label{thm:fourier}

\begin{enumerate} 

\item{For any $a \in \Z$, $ \F^2(a) = S^{-1}(a)$; }
\item{Let $\eta$ denote the isomorphism of commutative algebras $\eta : C_0 
\to C_r$,
$\eta(f) = f \leftharpoonup K_{2\rho}$, where $f(\leftharpoonup K_{2\rho})
(x) = 
f(K_{2\rho}x)$ for any $x \in \Uf$. Then for any $f \in C_0$, we have
$ \eta^{-1} \circ {\F'}^2 \circ \eta (f) = f \circ S $.}
\end{enumerate}
\end{Thm}

\begin{proof}
(1) is equivalent to the analogous statement in \cite{LM}: \\
Set $\mathfrak{S}_- (a) \equiv m({\id} \otimes \lambda_l)(R_{12}^{-1}(1 
\otimes a) R_{21}^{-1})$,
$a \in \Z$, then ${\mathfrak{S}_-}^2 = S^{-1}$. 

Let $R = \sum a_i \otimes b_i$. For a central element $a$ write 
$$ \F(a) = m(\lambda_r \otimes \id)(S(a) \otimes 1)(R_{21}R_{12}) = $$
$$ = m(\lambda_r \otimes \id)(S \otimes \id)((S^{-1} \otimes \id)(\sum b_i 
a_j \otimes a_i b_j)
(a \otimes 1)) = $$
$$ = m(\lambda_r \circ S \otimes \id)((\sum S^{-1}(a_j) S^{-1}(b_i) \otimes 
a_i b_j)(a \otimes 1)) = $$
$$ = m(\lambda _l \otimes \id)((\sum S^{-1}(b_i)S(a_j) \otimes a_i b_j)(a 
\otimes 1)) =$$
$$ = m(\lambda_l \otimes \id)((R_{21}^{-1} R_{12}^{-1})(a \otimes 1)) =
m(\id \otimes \lambda_l)((R_{12}^{-1} R_{21}^{-1})(1 \otimes a)) = 
{\mathfrak S}_-(a). $$

Here we have used the fact that $\lambda_r \circ S$ is a left integral, the 
property of 
left-invariant functionals $\lambda_l(xy) = \lambda_l(y S^2(x))$ and the 
identities 
$(S \otimes \id)(R) = (\id \otimes S^{-1})(R) = R^{-1}$ \cite{Dr}. 

(2) follows from (1). For $f_r \in C_r$ write 
$$ {\F'}^2 (f_r) = \J^{-1} \circ {\F}^2 \circ \J (f_r) = 
\J^{-1} \circ S^{-1} \circ \J (f_r) = $$
$$ = \J^{-1} \{ m(f_r \otimes S^{-1})(R_{21}R_{12}) \}= 
\J^{-1}\{ m(f_r \circ S \otimes \id)(S^{-1} \otimes S^{-1})(R_{21} R_{12}) 
\} = $$
$$ = \J^{-1}\{m(f_r \circ S \otimes \id)(R_{12} R_{21})\} = \J^{-1}\{m(\id 
\otimes f_r \circ S)
(R_{21} R_{12}) \}, $$

where we have used $(S^{-1} \otimes S^{-1})(R) =R$. For $f_r \in C_r$, 
$f_r \circ S \in C_l$.
Consider the algebra isomorphism from $C_l$ to $C_r$, given by the action 
of  
$K_{4\rho}$. Then by \cite{Dr} for any $g \in C_l$,
$(\id \otimes g)(R_{21} R_{12}) = ((g \leftharpoonup K_{4\rho}) \otimes 
\id)(R_{21} R_{12})$. 
Therefore we have 
$$ {\F'}^2 (f) = \J^{-1}\{m((f_r \circ S) \leftharpoonup K_{4\rho}) \otimes 
\id)(R_{21} R_{12}) =
(f_r \circ S) \leftharpoonup K_{4\rho}.$$

Now let $f_r = \eta (f)$, $f \in C_0$. Then  
$$ \eta^{-1} \circ {\F'}^2 \circ \eta (f) = 
((f \leftharpoonup K_{2\rho})\circ S) \leftharpoonup K_{2\rho} 
 = f \circ S, $$
as required. 
\end{proof}
In particular, \thm{fourier} implies that 
the square of the Fourier transform, conjugated by $\eta$, maps 
the character of a simple $\Uf$-module to the character of its dual: 
$$(\eta^{-1} \circ {\F'}^2 \circ \eta) ({\rm ch}L(\mu)) = {\rm ch}(L(\mu)
^*).$$

\begin{Thm} \label{thm:intersec}
$$ \tZ \cap \pZ = {\rm Ann}\,({\rm Rad}\;\tZ) . $$ 
\end{Thm}
  
The proof is based on the following two lemmas

\begin{Lem} \label{lem:St}
 The character of the Steinberg module 
$St = L((l-1)\rho)$ annihilates the radical of $\R$.
\end{Lem}

\begin{proof}
Consider the block decomposition of the center $\Z = \oplus_{\mu \in \tmX} 
\Z_\mu$.
The block $\Z_{(l-1)\rho} =\Z_{St}$ consists of 
elements which act nontrivially only on indecomposable modules with highest 
weights in 
the $\bullet$-orbit of $(l-1)\rho$. Since this weight is stabilized by the 
$\bullet$ action, 
and the Steinberg module is simple, projective and injective in the category 
of finite dimensional modules over $\Uf$ \cite{APW2}, the block $\Z_{St}$ 
is one-dimensional. By the structure  
of $\pZ$ (\prop{pZ}) and $\tZ$ (\thm{BG}), 
we have $\tZ_{St} = \pZ_{St} = \Z_{St}$. Let $z_{St}=\phi^{-1}(ch_{q^{-1}}St) 
\in \Z_{St}$. Then 
$z_{St} \in \tZ \cap \pZ$.

Write the sequence of mappings 
$$\R_r \stackrel{\J}{\longrightarrow} \tZ \stackrel{\phi}{\longrightarrow} 
\phi(\tZ) 
\stackrel{\J}{\longrightarrow} \J \circ \phi (\tZ) \stackrel{\phi}
{\longrightarrow} 
(\phi \circ \J)^2(\R_r) = \R_r. $$ 

The first arrow is the definition of $\J$, and the last equality follows by 
\thm{fourier},
since $\F' = \phi \circ \J$. Then we have $\phi^{-1}({\F'}^2(\R_r))=\phi^{-
1}(\R_r)=\pZ$ 
by the definition of $\phi$, and hence $\J \circ \phi (\tZ) = \F (\tZ) = \pZ$.
Also, $\F(\pZ) = \J \circ \phi (\pZ) = \J (\R_r) = \tZ$. 

Now since $z_{St} \in \tZ \cap \pZ$, we have 
$\F (z_{St}) \in \F(\tZ) \cap \F(\pZ) = \pZ \cap \tZ$. 
The ideal $\pZ$ annihilates the radical of $\Z$ by \prop{pZ}, and therefore 
its intersection with 
$\tZ$ belongs to the annihilator of the radical of $\tZ$. 
We obtain $\F(z_{St}) = \J \circ \phi (z_{St}) \in {\rm Ann}({\rm Rad}\;
\tZ)$. 
The map $\J$ is an algebra isomorphism between $\R_r$ and $\tZ$, and 
therefore 
$\phi (z_{St}) \in {\rm Ann}\,({\rm Rad} \R_r)$. By the definition of 
$z_{St}$, $ch_{q^{-1}}St = \phi(z_{St})$, and hence this element 
annihilates the radical of $\R_r$, or 
equivalently,
$ch St$ annihilates the radical of $\R$. Here we used the isomorphism of 
algebras 
$\R \to \R_r$, given by the action of $K_{-2\rho}$ (\sec{Hopf}). 
\end{proof}

\begin{Lem} \label{lem:St^2}
 The square of the Steinberg character 
spans the square of the annihilator of the radical of $\R$.
\end{Lem}
\begin{proof}
 Recall that $\R$ is a direct sum of local Frobenius algebras. Its socle 
(annihilator of the radical) is spanned by the 
socles of each of the $|\tmX|$ blocks. Only one of these blocks is 
semisimple, the one corresponding to the orbit of the Steinberg weight
$(l-1)\rho$. The socle of this block is an idempotent, while 
the socles of all other blocks are nilpotent of second degree. Therefore, 
the 
square of the socle of $\R$ is one-dimensional. Since 
${\mathrm ch}St \in {\rm Soc}(\R)$ by \lem{St}, 
$({\mathrm ch}St)^2 \in ({\rm Soc}(\R))^2$ and it is nonzero as 
a character of a direct sum of modules. Therefore, 
$({\mathrm ch}St)^2$ spans $({\rm Soc}(\R))^2$. 
\end{proof}

\begin{proof} {\it of the Theorem.}

Obviously $\tZ \cap \pZ \subset {\rm Ann}\,({\rm Rad }\;\tZ) \equiv {\rm 
Soc}\,\tZ$. We will show that 
${\rm Ann}\,({\rm Rad }\;\tZ) \subset \pZ$.

For each $\mu \in \tmX$, denote by $z_\mu$ an element of the one-
dimensional 
socle of the block $\tZ_\mu$. We want to show that 
each $z_\mu$ belongs to $\pZ$, or equivalently that $\phi(z_\mu)$ is in 
$\R_r$. 
We will use the lemmas above, recalling that $\R \simeq \R_r$, the 
isomorphism given by the action of $K_{2\rho}$ which maps characters of 
$\Uf$-modules to $q^{-1}$-characters.  

As before let $z_{(l-1)\rho} \equiv z_{St} = \phi^{-1} (ch_{q^{-1}}St)$. 
Then the element $z_{St}$ spans 
the square of ${\rm Soc} \,\tZ$, and therefore, by \lem{St^2} the isomorphism 
$\J$ maps $ch_{q^{-1}} St^2$ to the subspace spanned by $z_{St}$.
By \lem{St}, ${ch}_{q^{-1}}St$ is mapped by $\J$ to an element 
$\sum_{\nu \in \tmX} a_\nu z_\nu \in {\rm Soc}\,\tZ$ for some coefficients
$a_\nu \in {\mathbb C}$. Then up to nonzero scalar multiples 
one can write the following sequence of  mappings: 
$$ \;\;\;\;\;\;\;\;\;\;\; {ch}_{q^{-1}}St^2 \stackrel\J{\longrightarrow} z_
{St} \stackrel{\phi}{\longrightarrow} 
 {ch}_{q^{-1}}St \stackrel\J{\longrightarrow}$$
$$  \stackrel\J{\longrightarrow} \sum_{\nu \in \tmX} a_\nu z_\nu
\stackrel{\phi}{\longrightarrow}
 \sum_{\nu \in \tmX} a_\nu p(\nu)  = {\mathrm ch}_{q^{-1}}St^2. $$

Here set $p(\nu) \equiv 
\phi(z_\nu)$. 
 The last 
equality is a consequence of \thm{fourier} since ${\mathrm ch}St$ 
and its square are invariant 
under the antipode. 
Since the square of ${\rm Soc}(\R_r)$ is one-dimensional, 
${ch}_{q^{-1}}St^2$ coincides up to a scalar multiple with the product 
${ch}_{q^{-1}}St^2 \cdot f$ for any $f \in \R_r$. Therefore by weight 
considerations, ${\mathrm ch}_{q^{-1}}St^2$ decomposes into a sum 
containing elements in all blocks of the category, and hence all $a_\nu$ in 
$\sum_{\nu \in \tmX} a_\nu p(\nu) = {\mathrm ch}_{q^{-1}}St^2$ are nonzero. 
 Act on both sides of the last equality  
by an idempotent of one block $1_\mu$. 
Then we get $a_\mu p_\mu = pr_\mu ({\mathrm ch}_{q^{-1}}St)^2$, the 
projection of the $q^{-1}$-character of $St^2$ to the block $(\R_r)_\mu$, 
which is spanned by the $q^{-1}$-characters of simple modules with highest 
weights in the $W\bullet$ orbit of $\mu$. This means that 
each $p(\mu)$ belongs to $\R_r$. 
 This completes the proof of the Theorem.   
\end{proof}

\begin{Rem}
The subalgebra $\tZ \cap \pZ$ is an ideal in $\Z$ closed under the quantum 
Fourier transform, which suggests that it has an interesting representation 
theoretical interpretation. In fact, one can show that it is isomorphic to 
the ideal of characters of the subcategory of projective modules over $\Uf$. 
The properties of $\tZ \cap \pZ$ are studied in detail in \cite{La}. 
 \end{Rem}

Now we combine the results of \sec{Hopf}, \sec{tZ} and \sec{inters} to 
formulate the main theorem on the structure of a subalgebra in the center 
of $\Uf$:

\begin{Thm} \label{thm:main}
The center $\Z$ of the small quantum group $\Uf$ contains a subalgebra 
$\tZ + \pZ,$ which is invariant under the quantum Fourier transform, of 
dimension 
$${\rm dim}(\tZ + \pZ) = \sum_{\mu \in \tmX}( 2[W:W_\mu]-1) =2l^r-|\tmX|,
$$  
where $\tmX \subset (P/lP)_+$ enumerates the blocks of the category of 
finite dimensional $\Uf$-modules. The intersection $\tZ \cap \pZ$ is an 
ideal in $\Z$ isomorphic to ${\rm Soc}\,\tZ$. It contains one element in 
each block of $\Z$, has dimension ${\rm dim}(\tZ \cap \pZ) = |\tmX|$, 
and is also invariant under the quantum Fourier transform. The ideal $\pZ 
\subset \Z$ is isomorphic to ${\rm Soc}\,\Z$. The subalgebra 
$\tZ \subset \Z$ is isomorphic to the algebra $\R$ of characters of finite 
dimensional $\Uf$-modules. The dimensions of $\tZ$ and $\pZ$ coincide and 
are equal to $\sum_{\mu \in \tmX}[W:W_\mu] = l^r$.
\end{Thm}

\begin{Rem} 
The constructed subalgebra $\tZ +\pZ$ provides a supply of central 
elements, while the mappings $\F$, $\phi$ can serve as tools to 
describe them. For example, it is clear that the 
two-sided integral $\Lambda$ of $\Uf$, which is not an element of $\tZ$,
is the Fourier transform of the unit, 
$$ \Lambda = \F(1). $$

Another possible application originates from the theory of quantum 
topological invariants 
of $3$-manifolds 
\cite{RT}, \cite{Hen} \cite{Lyu}, \cite{Ker2}.  
It is known that a quantum group at a root of unity gives rise to a family 
of quantum 
surgical invariants, computed by assigning an element of the algebra to 
each component 
of the knot, using the quasitriangular structure to unknot the 
interlacings, and then 
evaluating the result against a functional in the dual quantum group. 
The invariant of the Reshetikhin-Turaev type is generated by a functional 
$\mu_{RT}$ which is a linear combination of $q^{-1}$-characters of 
certain simple modules \cite{RT}. At the same time, the right integral 
$\lambda_r$ of a quazitriangular finite dimensional Hopf algebra yields the 
Kauffman-Radford invariant \cite{KR}. 
It was pointed out in 
\cite{Hen} that all surgical invariants associated to a finite dimensional 
quantum group can be 
obtained from the Kauffman-Radford invariant by the 
left action of certain central elements, in particular  
$\mu_{RT} = z_{RT}
\rightharpoondown \lambda_r$.    
In our case it means that  $z_{RT} = \phi^{-1}(\mu_{RT})$, where 
$\mu_{RT} \in \R_r$. 
Then \prop{pZ} 
implies that $z_{RT} \in \pZ = 
{\rm Soc}\Z$ and $z_{RT}$ annihilates the 
radical of $\Z$ (cf. in \cite{Ker2}: $z_{RT}^2=0$). 
A more careful consideration based on the properties of $\mu_{RT}$ 
shows that $z_{RT}$ 
does not belong to the obvious subalgebra  $\tZ$ of the center. 
Although we consider here the case $\Uf$, the above arguments 
can be carried out in a more general situation, e.g. for a finite dimensional 
quantum group specialized at an even root of unity, which was the 
original setting for the Reshetikhin-Turaev invariant. One can deduce that
the two invariants are always related by the left action of an element 
which is not contained in the obvious part of the center 
and annihilates its radical.  

\end{Rem}

\section{The action of the center in projective modules} \label{sec:act}

In this section we derive a consequence of \thm{intersec}, which describes 
the action of certain central elements by the endomorphisms of projective 
modules. We will use freely the known facts from the representation theory 
of $\Uq$ and $\Uf$ in the proof below; a detailed exposition can be found 
in \cite{Lus2},
\cite{And}, \cite{APW1}, \cite{APW2}.

\begin{Thm} \label{thm:act}
 For an indecomposable projective $\Uf$- module $P((l-1)\rho +\mu)$, which 
is a restriction of the projective $\Uq$-module of highest weight 
$(l-1)\rho +\mu$ with  $\mu \in \tilde{X} = \{ \nu \in P_+ : \langle \nu, 
\alpha \rangle < l \; {\rm for} \; {\rm all} \; \alpha \in R_+ \}$, 
there exists a unique surjective homomorphism from $\tZ$ to 
${\rm End}_{U_q^{fin}} P((l-1)\rho +\mu)$, which is an isomorphism on the 
block 
$\tZ_\mu$, and zero on all other blocks.
\end{Thm}

\begin{proof}
 The restriction of $P((l-1)\rho +\mu)$ to $\Uf$ is a projective 
indecomposable module and the projective cover of $L((l-1)\rho +w_0(\mu))$ 
\cite{APW2}. 
By the definition of the block $\tZ_\mu$, its elements act as endomorphisms 
of a projective module with composition factors in the corresponding 
$W\bullet$-orbit of $\mu$,
while elements of the other blocks ${\Z}_{\nu \neq \mu}$ act on it by zero. 
Let $z_\mu \in {\rm Soc}\,\tZ_\mu$. 
By \lem{St}, $ch St \cdot \R \subset {\rm Ann}({\rm Rad}\;\R)$. For any 
$\nu \in (P/lP)_+$, it is possible to find $f \in \R$ such that $ch St 
\cdot f$ contains $ch L(\nu)$ with a nonzero coefficient in its 
decomposition with respect to the basis of simple characters. Therefore, 
any character of a simple module appears in some combination in 
the linear span of ${\rm Ann}
({\rm Rad}\;\R)$. 
This means that $z_\mu$ acts nontrivially in each projective cover 
of a simple module with the highest weight in the $W\bullet$-orbit of 
$\mu$, mapping it to its simple socle.
 
Since $\tZ$ is Frobenius, for any $z_1 \in  \tZ_\mu$ there exists 
$z_2 \in  \tZ_\mu$ such that $z_1 z_2 = z_\mu$. Then the kernel of the 
mapping 
of $\tZ_\mu$ to ${\rm End}_{\Uf} P((l-1)\rho +\mu)$ is trivial. 
Now compare the dimensions. 

First note that ${\rm End}_{U_q^{res}} P((l-1)\rho +\mu)= [W:W_\mu]$, where 
$W_\mu$ is the 
stabilizer of $\mu$ with respect to the $W\bullet$ action. This follows 
from the structure of filtrations of projective modules over $\Uq$. Each 
projective module has a filtration by Weyl modules $W(\nu)$, which are the 
universal highest weight modules for $\Uq$ \cite{APW1}. Then we have the 
reciprocity relation for the projective cover of a simple module 
$[ P((l-1)\rho +\mu) : W(w\cdot \mu) ]=[ W(w \cdot \mu) : L((l-1)\rho +w_0
(\mu)) ]$ (\cite{AJS}, \S 4.15). 
The Weyl filtration for a projective module in 
the $\tilde{X} +(l-1)\rho$ alcove has the following form:
 $$ {\mathrm ch}P((l-1)\rho +\mu) = \sum_{w \in W/W_\mu} {\mathrm ch}W((l-
1)\rho +w(\mu)).$$ 
The last formula of course can be deduced from the general result on the 
characters of 
tilting modules in \cite{Soe2} (since all projective modules over $\Uq$ are 
tilting, \cite{APW2}). It also follows by 
direct computation, using the fact that any projective module appears as a 
direct summand of a  
tensor product $St \otimes E$ for some finite dimensional module $E$ (\cite
{APW1}).

We need to check that there are no extra endomorphisms of the restriction 
of $P((l-1)\rho +\mu)$ over $\Uf$. Suppose there are additional composition 
factors in $P((l-1)\rho +\mu)$ isomorphic to its maximal semisimple quotient 
$L((l-1)\rho +w_0(\mu))$. Then by the tensor product decomposition 
for simple modules (\cite{Lus1}), a simple $\Uq$-module $L((l-1)\rho +w_0
(\mu) +l\nu)$ has to appear in the filtration of $P((l-1)\rho +\mu)$ for 
some $\nu \in P_+$. By the strong linkage principle (Thm. 8.1, \cite
{APW1}), this means
that the weight $ (l-1)\rho +w_0(\mu) +l\nu$ should be less than or equal to  
at least one of the weights $\{(l-1)\rho +w(\mu)\}_{w \in W}$ of the Weyl 
composition factors in $P((l-1)\rho +\mu)$, and belong to the same $\Wl 
\cdot$-orbit. For $\mu = \rho$, we should 
have then 
$$ (l-1)\rho +w(\rho) -((l-1)\rho -\rho +l\nu)=w(\rho)+\rho-l\nu \in P_+,$$
which is false for any $w \in W$ and any dominant $\nu \neq 0$. Therefore, 
there are no composition factors of the form $L((l-1)\rho +w_0(\rho))$ in 
any of the $W((l-1)\rho +w(\rho))$. 
The statement for other weights $\mu \in \tilde{X}$ such that $(l-1)\rho +
\mu$ is regular, is derived using the translation principle \cite{APW1} 
which states that the structure of the filtration of a Weyl module is the 
same for all highest weights in the same alcove. The proof for the lowest 
walls of $(l-1)\rho + \tilde{X}$ is obtained similarly, replacing $\rho$ 
by a suitable fundamental weight $\omega_i$ in the argument above.  
This means that $ {\rm dim}({\rm End}_{U_q^{fin}} P((l-1)\rho +\mu)) = 
{\rm dim}({\rm End}_{U_q^{res}} P((l-1)\rho +\mu))$ for any $\mu \in 
\tilde{X}$.

Therefore, we get 
$$ {\rm dim}({\rm End}_{U_q^{fin}} P((l-1)\rho +\mu)) = [W:W_\mu] = {\rm 
dim}{\tZ}_\mu, $$ 
and the two algebras are isomorphic.
\end{proof} 

\begin{Cor} \label{cor:sl_2}
For $\Lg = {\Sl}_2$, the subalgebra $\tZ + \pZ$ 
coincides with the whole center of $\Uf$. 
\end{Cor} 

\begin{proof} We know that the subalgebra 
$\tZ$ contains all idempotents of $\Z$ (\prop{idemp}).
We will estimate from above the dimension of the radical of the center, 
computing  
the number of inequivalent nilpotent edomorphisms of all projective 
modules.

Fix $\mu \in (P/lP)_+$. If $\mu = (l-1)\rho$, then it is fixed by the 
$W\bullet$ action and the corresponding block contains the single 
projective module $St$ which is simple and has no nilpotent 
endomorphisms. Suppose that $\mu$ is regular, then its $W\bullet$ orbit 
contains two weights. Denote the corresponding indecomposable 
projective modules by $P_1$ and $P_2$. 
In the case $\Lg ={\Sl}_2$, all indecomposable projective modules have 
highest 
weights 
in $(l-1)\rho + \tilde{X}$. Therefore, they all satisfy the condition of 
\thm{act}.  Then since ${\rm dim}\tZ_\mu =2$ for any regular $\mu$ by \thm
{BG},  we have ${\rm dim}({\rm End} P_1) = {\rm dim}({\rm End} P_2)=2$, and 
each algebra of endomorphisms contains a nilpotent 
element and an identity. Therefore, the center of ${\rm End}(P_1\oplus P_2)$ 
contains at most two nilpotent elements.  
On the other hand, the dimension of the regular block $\pZ_\mu$ of 
$\pZ$ is equal to 
$|W| =2$, 
which means that all nilpotent central elements are contained in $\pZ$, and 
hence ${\Z} = {\tZ} + \pZ$. 
\end{proof}

\begin{Rem}
The above argument fails in higher rank. By (\cite{AJS}, \S 19)
 the restriction to $\Uf$ of a projective 
$\Uq$-module with a regular highest weight outside of $(l-1)\rho +\tilde{X}$, 
has more than $|W|$ linearly independent endomorphisms.
\end{Rem}

\section{Examples.} \label{sec:ex}
\subsection{The case $ \Lg = {\Sl}_2$.}

The center $\Z$ decomposes as a sum of ideals parametrized by the set
 $\tmX = \{0, \ldots (\frac{l-3}{2}), (l-1)\}$. A regular orbit contains 
two weights 
$\{ i, (l-2-i) \}_{i =0}^{\frac{l-3}{2}}$. The only singular orbit contains 
the 
Steinberg weight $(l-1)$. Therefore, 
$${\Z} = \oplus_{i=0}^{\frac{l-3}{2}} {\Z}_i \oplus {\Z}_{l-1}. $$
The subalgebra $\tZ = \J(\R_r) \subset \Z$,
is isomorphic to the algebra $\R$,
$$\R \simeq {\mathbb C}[x]^W \otimes_{{\mathbb C}[x^l]^W}{\mathbb C} \simeq 
{\mathbb C}[x + x^{-1}]/\langle x^l + x^{-l} -2 \rangle, $$ 
which is spanned by the characters of simple modules ${\mathrm ch}L(i) 
\equiv \xi(i)$ for 
$i = 0, \ldots ,(l-1)$.
 
 \thm{BG} specializes to 
$$ \tZ = \oplus_{i=0}^{\frac{l-3}{2}} \tZ_i 
\oplus \tZ_{l-1} = \oplus_{i =0}^{\frac{l-3}{2}} 
{\mathbb C}[x] \otimes_{{\mathbb C}[x + x^{-1}]} {\mathbb C} \oplus 
{\mathbb C}. $$

This means that each regular block ${\Z}_i$ contains a two-dimensional 
subalgebra isomorphic to 
${\mathbb C}[(x-1)]/\langle (x-1)^2 \rangle$, and ${\Z}_{l-1} = \tZ_{l-1} = 
{\mathbb C}$. 

By \prop{pZ}, the socle of $\Z$ coincides with the inverse image of the 
Grothendieck ring under the isomorphism $\phi$, which preserves the block 
structure.
Therefore, each regular ${\Z}_i$ contains two linearly independent elements 
of the 
socle, corresponding to the two simple characters in the orbit: $\phi^{-1}
(\xi(i))$ and $\phi^{-1}(\xi(l-2-i))$; their product and 
both squares are zero, since the socle of a regular block contains no 
idempotents. 
 The singular block ${\Z}_{l-1}$ consists of one idempotent element 
$\phi^{-1}(\xi(l-1))$. 

By \thm{BG}, the subspace ${\rm Soc}\,\tZ_i$ is one-dimensional for each 
block. By \lem{St}, its image in $\R$ contains the characters of tensor 
products of the Steinberg module with finite dimensional modules. For $\Lg = 
{\Sl}_2$, $ch St \cdot \R$ is spanned by the characters of the tilting 
modules
${\mathrm ch}T(2l-2-i) = {\mathrm ch}T(l+i)=2(\xi(i) + \xi(l-2-i))$, and 
therefore
${\rm Soc}\,\R$ is spanned by 
$$\{ (\xi(i) + \xi(l-2-i)), \xi(l-1) \}_{i=0}^{\frac{l-3}{2}}.$$ 
Its dimension is $\frac{l+1}{2}$. 

By \cor{sl_2}, the whole center $\Z$ coincides with the sum of the two 
subalgebras $\tZ + \pZ$ in this case. 
Therefore, we have 
$$ \Z \simeq \tZ+ \pZ \simeq \oplus_{i =0 \ldots \frac{l-3}{2}, i= l-1} 
(\tZ+ \pZ)_i,$$
where $(\tZ+ \pZ)_i \simeq \{1_i, t_i, s_i \}$ with 
$1_i$ acting as the unit of the block, and $t_i s_i = t_i^2 = s_i^2 =0$ for 
$i =0, \ldots 
\frac{l-3}{2}$, and $(\tZ+ \pZ)_{l-1} \simeq \{1_{l-1} \}$. 
This coincides with the description of the center obtained in \cite{Ker1} 
in this particular case.

\subsection{The case $\Lg = {\Sl}_3$.}

The set 
$\tmX$ consists of a one-dimensional orbit corresponding to the Steinberg 
module; $(l-1)$ $3$-dimensional orbits corresponding to the weights 
stabilized by one reflection in $W\bullet$; and $\frac{(l-2)(l-1)}{6}$  
$6$-dimensional orbits corresponding to $W\bullet$-regular weights. 

Check the dimension: 
$$1 + 3(l-1) + 6 \frac{(l-2)(l-1)}{6} = l^2 = {\rm dim}\R.$$ 

\thm{BG} specializes to 
$$ \tZ \simeq \bigoplus_{i=1}^{\frac{(l-1)(l-2)}{6}} 
{\mathbb C}[x,y,z]/I_1 \oplus \bigoplus_{i=1}^{l-1} {\mathbb C}[x]/I_2 
\oplus 
{\mathbb C}, $$
where $I_1 = \langle x+y+z, xy+yz+xz, xyz \rangle$, and 
$I_2 = \langle x^3 \rangle$. 

The algebra $\pZ$ has the same block decomposition as $\tZ$ and the same 
dimension of each block. Denote by 
$C_n$ the $n$-dimensional vector space with zero multiplication.  Then 
$$\pZ \simeq  \oplus_{i=1}^{\frac{(l-1)(l-2)}{6}}C_6 \oplus 
\oplus_{i=1}^{l-1} C_3 \oplus {\mathbb C}. $$ 
The multiplication between $\tZ$ and $\pZ$ is determined by the condition 
that $\pZ \simeq {\rm Soc}\,\Z$. 

The dimension of $\tZ + \pZ$ is equal to 
$$ 1 + 5(l-1) + 11\frac{(l-1)(l-2)}{6}, $$ 
which is an integer since $(l, 3)=1$ by the assumption on $l$.


\end{document}